\documentclass[10pt]{article}
\usepackage{amssymb}
\usepackage{amsmath}
\usepackage[dvips]{graphics}
\usepackage{epsfig}

\begin{document}
\title{On the order and the type of an entire function}

\author{\textbf{Vassilis G. Papanicolaou$^1$, Eva Kallitsi$^{\,2}$, and George Smyrlis$^3$}
\\\\
Department of Mathematics
\\
National Technical University of Athens,
\\
Zografou Campus, 157 80, Athens, GREECE
\\
{\tt $^1$papanico@math.ntua.gr, $^2$evpapa@hotmail.com,}
\\
{\tt $^3$gsmyrlis@math.ntua.gr}}
\maketitle

\begin{abstract}

In this short article we present some properties regarding the order and the type of an entire function.
\end{abstract}

\\\\
\textbf{2010 AMS Mathematics Classification.} 30D20.

\section{Introduction}
Let
\begin{equation}
g(z) = \sum_{n=0}^{\infty} a_n z^n,
\qquad
z \in \mathbb{C},
\label{PREL1}
\end{equation}
be an entire function and
\begin{equation}
M(r) = M_g(r) := \sup_{|z| \leq r} |g(z)| = \max_{|z| = r} |g(z)|,
\qquad
r > 0,
\label{PREL2}
\end{equation}
its maximum modulus.

Recall that the order of $g(z)$ is the quantity \cite{H}
\begin{equation}
\rho = \rho(g) := \limsup_{r \to \infty} \frac{\ln\ln M(r)}{\ln r}.
\label{PREL3}
\end{equation}
In other words, the order $\rho$ of $g(z)$ is the smallest exponent $\rho' \geq 0$ such that for any given $\varepsilon > 0$ there is a
$r_0 = r_0(\varepsilon) > 0$ for which
\begin{equation}
|g(z)| \leq \exp\left(|z|^{\rho' + \varepsilon}\right)
\qquad \text{whenever }\;
|z| \geq r_0.
\label{PREL4}
\end{equation}
Clearly, $0 \leq \rho \leq \infty$.

Let us also recall \cite{H} that if $0 < \rho < \infty$, the quantity
\begin{equation}
\tau = \tau(g) := \limsup_{r \to \infty} \frac{\ln M(r)}{r^{\rho}}
\label{PREL5}
\end{equation}
is the type (of the order) of $g(z)$. In other words, $\tau$ is the smallest number $\tau' \geq 0$ such that for any given
$\varepsilon > 0$ there is a $r_0 = r_0(\varepsilon) > 0$ for which
\begin{equation}
|g(z)| \leq \exp\Big((\tau' + \varepsilon) |z|^{\rho}\Big)
\qquad \text{whenever }\;
|z| \geq r_0.
\label{PREL6}
\end{equation}
Clearly, $0 \leq \tau \leq \infty$. If $\tau = 0$, we say that $g(z)$ is of \emph{minimal type}, whereas if $\tau = \infty$, we say that $g(z)$ is of
\emph{maximal type}. In the extreme cases where $\rho = 0$ or $\rho = \infty$ the type is not defined.

It follows easily from \eqref{PREL4} and \eqref{PREL6} that for any fixed $\zeta \in \mathbb{C}$ the entire functions
\begin{equation}
g(z)
\quad \text{and} \quad
g_{\zeta}(z) := g(z + \zeta)
\qquad\quad \text{have the same order and type}.
\label{PREL7}
\end{equation}

A well-known fact of complex analysis is \cite{H} that the order $\rho$ and the type $\tau$ of $g(z)$ are given by the formulas
\begin{equation}
\rho = \limsup_n \frac{n \ln n}{-\ln |a_n|}
\label{PREL8}
\end{equation}
and (in the case where $0 < \rho < \infty$)
\begin{equation}
\tau = \frac{1}{e \rho}\limsup_n n |a_n|^{\rho/n}
\label{PREL9}
\end{equation}
respectively, where $a_n$, $n = 0, 1, \ldots$, are the coefficients of the power series of $g(z)$ as seen in \eqref{PREL1}.

%
%

\medskip

\textbf{Remark 1.} (i) Suppose $0 \leq \rho < \infty$. Then \eqref{PREL8} implies that for all $n$ sufficiently large we have
\begin{equation}
\frac{n \ln n}{-\ln |a_n|} \leq \rho + \varepsilon_n
\label{A1}
\end{equation}
for some sequence $\varepsilon_n$ of positive numbers with $\varepsilon_n \to 0$.

Formula \eqref{A1} can be written as
\begin{equation}
\ln\left(|a_n|^{1/n}\right) \leq -\frac{\ln n}{\rho + \varepsilon_n}
\label{A2}
\end{equation}
and for $\rho' > \rho$ \eqref{A2} yields
\begin{equation}
\ln\left(n \, |a_n|^{\rho'/n}\right) \leq -\left(\frac{\rho'}{\rho + \varepsilon_n} - 1\right) \ln n,
\label{A3}
\end{equation}
which implies (since, eventually, $\rho' > \rho + \varepsilon_n$)
\begin{equation}
\lim_n n |a_n|^{\rho'/n} = 0
\qquad \text{for }\;
\rho' > \rho.
\label{A4}
\end{equation}

(ii) Suppose $0 < \rho \leq \infty$. Then from \eqref{PREL8} we have that there is a subsequence $a_{n_k}$ such that
\begin{equation}
\lim_k \frac{n_k \ln n_k}{-\ln |a_{n_k}|} = \rho.
\label{A5}
\end{equation}
It follows that
\begin{equation}
\frac{\ln |a_{n_k}|}{n_k \ln n_k} = -\frac{1}{\rho} + o(1)
\quad \Longleftrightarrow \quad
\ln\left(|a_{n_k}|^{1/n_k}\right) = -\frac{1}{\rho} \ln n_k + o(\,\ln n_k).
\label{A6}
\end{equation}
Thus, for $0 < \rho' < \rho$ formula \eqref{A6} yields
\begin{equation}
\ln\left(n_k \, |a_{n_k}|^{\rho'/n_k}\right) = \left(1 - \frac{\rho'}{\rho}\right) \ln n + o(\,\ln n_k).
\label{A7}
\end{equation}
Therefore, $\lim_k n_k \, |a_{n_k}|^{\rho'/n_k} = \infty$ and, consequently,
\begin{equation}
\limsup_n n |a_n|^{\rho'/n} = \infty
\qquad \text{for }\;
\rho' \in (0, \rho).
\label{A8}
\end{equation}
Formulas \eqref{A4} and \eqref{A8} should be compared with formula \eqref{PREL9}.

\medskip

\textbf{Remark 2.} Let $g(z)$ of \eqref{PREL1} be entire of order $\rho \in (0, \infty)$, so that its type $\tau$ is defined, with
$0 \leq \tau \leq \infty$.

Suppose $0 <\tau < \infty$ and let $a_{n_k}$ be a subsequence for which the $\limsup$ in \eqref{PREL9} is attained, i.e.
\begin{equation}
\tau = \frac{1}{e \rho}\lim_k n_k |a_{n_k}|^{\rho/n_k}.
\label{A9}
\end{equation}
Then,
\begin{equation}
|a_{n_k}|^{\rho/n_k} = \frac{e \rho \tau}{n_k} + o\left(\frac{1}{n_k}\right) = \frac{e \rho \tau}{n_k} \big[1 + o(1)\big].
\label{A10}
\end{equation}
Taking logarithms in \eqref{A10} yields
\begin{equation}
\frac{\rho}{n_k} \ln |a_{n_k}| = \ln(e \rho \tau) - \ln(n_k) + o(1)
\label{A11}
\end{equation}
or
\begin{equation}
\frac{\rho \ln |a_{n_k}|}{n_k \ln n_k} = -1 + O\left(\frac{1}{\ln n_k }\right).
\label{A12}
\end{equation}
Therefore,
\begin{equation}
\lim_k \frac{n_k \ln n_k}{-\ln |a_{n_k}|}  = \rho,
\label{A13}
\end{equation}
i.e. the $\limsup$ in \eqref{PREL8} is also attained for the same subsequence $a_{n_k}$.

Notice that the converse is not true in general. If $a_{n_k}$ is a subsequence for which the $\limsup$ in \eqref{PREL8} is attained, then
the $\limsup$ in \eqref{PREL9} may not be attained for $a_{n_k}$ (e.g., take $g(z) = \sin z + \cos(2z)$ and $a_{n_k} = a_{2k+1}$).

If $\tau = \infty$ and $a_{n_k}$ is a subsequence for which the $\limsup$ in \eqref{PREL9} is attained, then for any $M > 0$ we have
\begin{equation}
|a_{n_k}|^{\rho/n_k} \geq \frac{M}{n_k}
\qquad \text{for all sufficiently large } k,
\label{A14}
\end{equation}
which implies
\begin{equation}
\frac{\rho \ln|a_{n_k}|}{n_k \ln n_k} \geq -1 + \frac{\ln M}{\ln n_k}
\label{A15}
\end{equation}
from which it follows that, again, the $\limsup$ in \eqref{PREL8} is also attained for the same subsequence $a_{n_k}$.

Finally, if $\tau = 0$, then \eqref{PREL9} implies that
\begin{equation}
\lim_n n |a_n|^{\rho/n} = 0.
\label{A16}
\end{equation}
In this case, though, there may exist subsequences $a_{n_k}$ for which the $\limsup$ in \eqref{PREL8} is not attained. For instance, let
$g(z) = g_e(z) + g_o(z)$, where $g_e(z)$ is an even entire function of order $\rho$ and type $0$, while $g_o(z)$ is an odd entire function of order
$< \rho$. Then $g(z)$ is of order $\rho$ and $0$ type, but the $\limsup$ in \eqref{PREL8} is not attained for the subsequence $a_{n_k} = a_{2k+1}$.

\medskip

Another tool that we will need in the sequel is the operator $(^{\sharp})$ defined as
\begin{equation}
g^{\sharp}(z) := \sum_{n=0}^{\infty} |a_n| z^n,
\label{PREL9a}
\end{equation}
when $g(z)$ is the entire function of \eqref{PREL1}. Notice that for any $r > 0$ we have
\begin{equation}
\max_{|z| \leq r} \left|g(z)\right| \leq g^{\sharp}(r) = \max_{|z| \leq r} \left|g^{\sharp}(z)\right|
\label{PREL9b}
\end{equation}
(the inequality can be strict).
Furthermore, since, in view of \eqref{PREL8} and \eqref{PREL9}, the order and type of $g(z)$ depend only on the sequence of the absolute values $\{|a_n|\}_{n \geq 0}$, they remain invariant under $(^{\sharp})$, i.e.
\begin{equation}
\rho\left(g^{\sharp}\right) = \rho(g)
\qquad \text{and} \qquad
\tau\left(g^{\sharp}\right) = \tau(g).
\label{PREL9c}
\end{equation}
Also,
\begin{equation}
\left(g'\right)^{\sharp}(z) = \sum_{n=0}^{\infty} |n a_n| z^{n-1} = \sum_{n=0}^{\infty} n |a_n| z^{n-1} = \left(g^{\sharp}\right)'(z),
\label{PREL9d}
\end{equation}
i.e. $(^{\sharp})$ commutes with the derivative operator.

Let us now set
\begin{equation}
a_n(z) := \frac{g^{(n)}(z)}{n!}
\qquad
n = 0, 1, \ldots
\label{PREL10}
\end{equation}
(so that $a_n(0) = a_n$). Then, in view of \eqref{PREL7}, formulas \eqref{PREL8}, \eqref{PREL9}, and \eqref{PREL10} yield
\begin{equation}
\rho = \limsup_n \frac{n \ln n}{-\ln |a_n(z)|}
\label{PREL11}
\end{equation}
and (in the case where $0 < \rho < \infty$)
\begin{equation}
\tau = \frac{1}{e \rho}\limsup_n n |a_n(z)|^{\rho/n}
= \frac{e^{\rho - 1}}{\rho}\limsup_n n^{1-\rho} \left|g^{(n)}(z)\right|^{\rho/n},
\label{PREL12}
\end{equation}
respectively, independently of the complex number $z$. An interesting question here is to inquire into the dependence on $z$ of the subsequence(s)
of $\{a_n(z)\}$ for which the $\limsup$ is attained in \eqref{PREL11} and \eqref{PREL12}.

Let us also notice that we can write \eqref{PREL11} in the equivalent form
(since $\lim_n |a_n(z)| = 0$ and, hence, $-\ln |a_n(z)|$ is eventually positive)
\begin{equation}
e^{-1/\rho} = \limsup_n |a_n(z)|^{\frac{1}{n \ln n}}
\label{PREL13}
\end{equation}
or, in view of \eqref{PREL10} and the fact that $\lim_n (n!)^{\frac{1}{n \ln n}} = e$,
\begin{equation}
\theta = \theta(\rho) := e^{1 - (1/\rho)} = \limsup_n \left|g^{(n)}(z)\right|^{\frac{1}{n \ln n}}.
\label{PREL14}
\end{equation}
Clearly, $\theta = \theta(\rho)$ is smooth and strictly increasing on $[0, +\infty]$, with $\theta(0) := \theta(0^+) = 0$ and
$\theta(+\infty) = e$.

Also, if
\begin{equation}
\theta^{\sharp} = \limsup_n \left|\left(g^{\sharp}\right)^{(n)}(z)\right|^{\frac{1}{n \ln n}},
\qquad
z \in \mathbb{C},
\label{PREL14a}
\end{equation}
then we must, obviously, have $\theta^{\sharp} = \theta$ since, as we have seen, $\rho\left(g^{\sharp}\right) = \rho(g) = \rho$. Thus, if we set
\begin{equation}
m_n(r) := \max_{|z| \leq r} \left|g^{(n)}(z)\right|,
\qquad
r > 0,
\label{PREL14b}
\end{equation}
then, clearly, $m_n(r) \leq (g^{(n)})^{\sharp}(r) = (g^{\sharp})^{(n)}(r)$ and, therefore
\begin{equation}
\theta = \limsup_n \left[m_n(r)\right]^{\frac{1}{n \ln n}} = \limsup_n \left[\left(g^{\sharp}\right)^{(n)}(r)\right]^{\frac{1}{n \ln n}},
\qquad
r > 0.
\label{PREL14c}
\end{equation}
Likewise, in view of \eqref{PREL12}, since $\tau(g^{\sharp}) = \tau(g) = \tau$ (in the case where $0 < \rho < \infty$) we have
\begin{equation}
\tau = \frac{e^{\rho - 1}}{\rho}\limsup_n n^{1-\rho} \left[m_n(r)\right]^{\rho/n}
= \frac{e^{\rho - 1}}{\rho}\limsup_n n^{1-\rho} \left[\left(g^{\sharp}\right)^{(n)}(r)\right]^{\rho/n},
\quad
r > 0.
\label{PREL14d}
\end{equation}

Now, let $\nu = \{n_k\}_{k=1}^{\infty}$ be a \emph{sequence of indices}, namely a strictly increasing sequence of positive integers. In the present work we are interested in the quantities
\begin{equation}
\rho_{\nu}(z) := \limsup_k \frac{n_k \ln(n_k)}{-\ln |a_{n_k}(z)|},
\qquad
z \in \mathbb{C},
\label{PREL15}
\end{equation}
and
\begin{equation}
\tau_{\nu}(z) := \frac{1}{e \rho}\limsup_k n_k \, |a_{n_k}(z)|^{\rho/n_k}
\qquad
z \in \mathbb{C}.
\label{PREL42}
\end{equation}
Let $\text{Ran}\,(\nu) := \{n_1, n_2, \ldots\}$ be the range of $\nu$ and $\mathbb{N} := \{1, 2, \ldots\}$ the set of natural numbers. If
$\mathbb{N} \setminus \text{Ran}\,(\nu)$ is a finite set, then it is clear from \eqref{PREL11} and \eqref{PREL12} that $\rho_{\nu}(z) = \rho$
and (if $0 < \rho < \infty$) $\tau_{\nu}(z) = \tau$ for all $z \in \mathbb{C}$. Therefore, to make things interesting we assume that
$\nu$ is a \emph{proper sequence of indices}, namely a sequence of indices such that the set $\mathbb{N} \setminus \text{Ran}\,(\nu)$ is infinite.
In this case it is obvious that there is a unique proper
sequence of indices $\mu=\{m_{\ell}\}_{\ell=1}^{\infty}$ such that $\text{Ran}\,(\mu) \cap \text{Ran}\,(\nu) = \emptyset$
and $\text{Ran}\,(\mu) \cup \text{Ran}\,(\nu) = \mathbb{N}$. We will call $\mu$ the \emph{complementary sequence of} $\nu$. It is clear that
for every $z \in \mathbb{C}$ we have
\begin{equation}
\rho = \max\{\rho_{\nu}(z), \, \rho_{\mu}(z)\}
\label{PREL17}
\end{equation}
and
\begin{equation}
\tau = \max\{\tau_{\nu}(z), \, \tau_{\mu}(z)\}.
\label{PREL44}
\end{equation}

Let us summarize the main results of the paper: Under the assumption that the sequence of indices $\nu=\{n_k\}_{k=1}^{\infty}$ satisfies the growth
condition $n_{k+1}/n_k \rightarrow 1$, in Section 2  we show that $\rho_{\nu}(z) = \rho$ for almost every $z \in \mathbb{C}$ and in Section 3 we show
that $\tau_{\nu}(z) = \tau$ for almost every $z \in \mathbb{C}$, provided that $\tau < \infty$; if $\tau = \infty$, then there is a dense $G_{\delta}$
subset $S$ of $\mathbb{C}$ such that $\tau_{\nu}(z) = \infty$ for $z \in S$.

These results, apart from being interesting per se, they can be used in the study of entire solutions of partial differential equations.

\section{Properties of the order}

For a proper sequence of indices $\nu=\{n_k\}_{k=1}^{\infty}$ we set (in the spirit of \eqref{PREL14})
\begin{equation}
\theta_{\nu}(z) := \exp\left(1 - \frac{1}{\rho_{\nu}(z)}\right) = \limsup_k \left|g^{(n_k)}(z)\right|^{\frac{1}{n_k \ln (n_k)}},
\qquad
z \in \mathbb{C}.
\label{PREL18}
\end{equation}
Thus, if $\mu$ is the complementary sequence of $\nu$, then \eqref{PREL17} gives
\begin{equation}
0 \leq \theta = \max\{\theta_{\nu}(z), \, \theta_{\mu}(z)\} \leq e.
\label{PREL20}
\end{equation}

We wish to determine how close is the quantity $\rho_{\nu}(z)$ to the order $\rho$ of $g(z)$ or, equivalently, how close is
the quantity $\theta_{\nu}(z)$ to the constant $\theta$ of \eqref{PREL14}.

Recall that a function $\phi(z)$, defined in a domain $\Omega$ of the complex plane and taking values in $\mathbb{R} \cup \{-\infty\}$, is called
subharmonic (in $\Omega$) if it is locally integrable and for any disk
$D_r(z_0) := \{z \in \mathbb{C} \, : \, |z - z_0| < r\} \subset \Omega$ we have
\begin{equation}
\phi(z_0) \leq \frac{1}{\pi r^2}\int_{D_r(z_0)} \phi(z)\, dx dy
\label{PREL21}
\end{equation}
(here, of course, $z = x + iy = (x, y)$ and the function $\phi(z)$ is subharmonic with respect to the real variables $x$ and $y$). Some authors
require \eqref{PREL21} to hold for almost every $z_0 \in \Omega$, in order to completely characterize subharmonic functions as functions whose
distributional Laplacian is nonnegative \cite{L-L}. However, such variants of the definition of subharmonicity are nonessential for our analysis.

If $A(z)$ is analytic in a domain $\Omega \subset \mathbb{C}$, then $\ln|A(z)|$ is subharmonic
in $\Omega$ (this follows, e.g., from the facts that (i) $\ln|z - z_0|$ is subharmonic and (ii) if $A(z)$ does not vanish in $\Omega$, then
$\ln|A(z)|$ is harmonic). Also, since $h(x) = e^{\alpha x}$, $x \in \mathbb{R}$, is convex for any $\alpha > 0$, Jensen's inequality implies that
$|A(z)|^{\alpha} = e^{\alpha\ln|A(z)|}$ too is subharmonic for any $\alpha > 0$.

\medskip

\textbf{Lemma 1.} The function $\theta_{\nu}(z)$ defined by \eqref{PREL18} is subharmonic in $\mathbb{C}$.

\smallskip

\textit{Proof}. Let us set
\begin{equation}
\Phi_n(z) := \sup_{k \geq n} \left|g^{(n_k)}(z)\right|^{\frac{1}{n_k \ln (n_k)}},
\qquad
n \geq 2.
\label{PREL22}
\end{equation}
Fix an $r > 0$ and restrict $z \in D_r := D_r(0) = \{z \, : \, |z| \leq r\}$. Then, $|g^{(n_k)}(z)| \leq (g^{\sharp})^{(n_k)}(r)$. Furthermore, as we have seen,
\begin{equation}
\theta = \limsup_n \left[(g^{\sharp})^{(n)}(r)\right]^{\frac{1}{n \ln n}} \leq e.
\label{PREL23}
\end{equation}
It follows that there is an $M = M(r) > 0$, such that $\Phi_n(z)$ of \eqref{PREL22} is $\leq M$ for all $n \geq 2$ and all $z \in D_r$.

Now, from the discussion preceding Lemma 1 we know that $|g^{(n_k}(z)|^{\frac{1}{n_k \ln (n_k)}}$ is subharmonic for any $k \geq 2$.
It, then, follows easily that $\Phi_n(z)$ is subharmonic in $D_r$ for every $n \geq 2$ (being finite and the supremum of a sequence of subharmonic
functions). Furthermore, it is obvious that $\Phi_n(z)$ decreases with $n$ and, in view of \eqref{PREL18},
\begin{equation}
\theta_{\nu}(z) = \lim_n \Phi_n(z),
\qquad
z \in D_r.
\label{PREL24}
\end{equation}
Therefore, by a simple application of the bounded convergence theorem we can conclude that $\theta_{\nu}(z)$ is subharmonic in $D_r$ and,
consequently, since $r$ is arbitrary, that $\theta_{\nu}(z)$ is subharmonic in $\mathbb{C}$.
\hfill $\blacksquare$

\medskip

\textbf{Remark 3.} It is a well-known fact \cite{L-L} that a subharmonic function $\phi(z)$ in a domain $\Omega$ is equal to an upper semicontinuous function
$\hat{\phi}(z)$ for almost every (a.e.) $z \in \Omega$. Therefore, Lemma 1 implies
\begin{equation}
\theta_{\nu}(z) = \hat{\theta}_{\nu}(z)
\qquad
\text{for \;a.e. } z \in \mathbb{C},
\label{PREL25}
\end{equation}
where $\hat{\theta}_{\nu}(z)$ is upper semicontinuous in $\mathbb{C}$.

\medskip

\textbf{Example 1.} Suppose $n_k = 2k$, $k = 1, 2, \ldots$, and $g(z) = \sin(\lambda z)$, where $\lambda \in \mathbb{C} \setminus \{0\}$. Then,
$\rho = 1$ and $\tau = |\lambda|$. Furthermore, since
\begin{equation*}
g^{(2k)}(z) = (-1)^k \lambda^{2k} \sin(\lambda z),
\end{equation*}
we have, in view of \eqref{PREL18},
\begin{equation*}
\theta_{\nu}(z)
= \left\{
  \begin{array}{cc}
    1, & \ \lambda z/\pi \in \mathbb{C} \setminus \mathbb{Z}; \\
    0, & \ \lambda z/\pi \in \mathbb{Z}, \\
  \end{array}
\right.
\end{equation*}
where $\mathbb{Z}$ is the set of integers. Obviously, $\theta_{\nu}(z)$ is subharmonic and it is equal to $\hat{\theta}_{\nu}(z) \equiv 1$ for all
except for countably many $z \in \mathbb{C}$.

\medskip

We are now ready for the main result of the section.

\medskip

\textbf{Theorem 1.} Let $\nu=\{n_k\}_{k=1}^{\infty}$ be a proper sequence of indices of \emph{subexponential growth}, namely
\begin{equation}
\frac{n_{k+1}}{n_k} \rightarrow 1
\qquad \text{as }\;
k \to \infty.
\label{PREL41a}
\end{equation}
Then, for an entire function $g(z)$ we have
\begin{equation}
\theta_{\nu}(z) = \theta
\qquad
\text{for \;a.e. } z \in \mathbb{C},
\label{PREL26}
\end{equation}
where $\theta$ and $\theta_{\nu}(z)$ are as in \eqref{PREL14} and \eqref{PREL18} respectively.

\smallskip

\textit{Proof}. By Lemma 1 we have that $\theta_{\nu}(z)$ is subharmonic in $\mathbb{C}$. Hence, in view of \eqref{PREL21} we must have
\begin{equation}
\theta_{\nu}(w) \leq \frac{1}{\pi r^2} \int_{D_r(w)} \theta_{\nu}(z) \, dx dy \leq \theta
\label{PREL27}
\end{equation}
for any $w \in \mathbb{C}$ and any $r > 0$.
Therefore, if for some $w \in \mathbb{C}$ we have that $\theta_{\nu}(w) = \theta$, then formula \eqref{PREL27} implies that $\theta_{\nu}(z) = \theta$
for a.e. $z \in D_r(w)$, which in turn implies $\theta_{\nu}(z) = \theta$ for a.e. $z \in \mathbb{C}$, since $r$ is arbitrary.

More generally, let us only assume that the supremum of $\theta_{\nu}(z)$ on some compact subset of $\mathbb{C}$ is $\theta$, namely that there is a
sequence $\{z_n\}_{n=1}^{\infty}$ with $\lim_n z_n = z_{\ast} \in \mathbb{C}$ and $\lim_n \theta_{\nu}(z_n) = \theta$. We will show that, we must, again,
have $\theta_{\nu}(z) = \theta$ for a.e. $z \in \mathbb{C}$.

Fix a disk $D_r(z_{\ast})$ and consider the disks $D_n := D_{r_n}(z_n)$, $n = 1, 2, \ldots$, so that $r_n$ is the largest radius satisfying
$D_n \subset D_r(z_{\ast})$. Using $w = z_n$ and $D_r(w) = D_n$ in \eqref{PREL27} yields
\begin{equation}
\theta_{\nu}(z_n) \leq \frac{1}{\pi r_n^2} \int_{D_n} \theta_{\nu}(z) \, dx dy
\leq \frac{1}{\pi r_n^2} \int_{D_r(z_{\ast})}  \theta_{\nu}(z) \, dx dy
\leq \frac{r^2}{r_n^2} \, \theta,
\qquad
n \geq 1,
\label{PREL28}
\end{equation}
thus, by letting $n \to \infty$ we obtain
\begin{equation}
\theta \leq \frac{1}{\pi r^2} \int_{D_r(z_{\ast})} \theta_{\nu}(z) \, dx dy
\leq \theta,
\label{PREL29}
\end{equation}
which tells us that $\theta_{\nu}(z) = \theta$ for a.e. $z \in D_r(z_{\ast})$ and, consequently, that $\theta_{\nu}(z) = \theta$ for a.e. $z \in \mathbb{C}$.


Finally, we will show that the assumption
\begin{equation}
\Theta_{\nu}(r) := \sup_{|z| \leq r} \theta_{\nu}(z) < \theta
\qquad \text{for every }\; r > 0
\label{PREL30}
\end{equation}
leads to a contradiction.

For a given $r > 0$ let us assume \eqref{PREL30} and fix an $\varepsilon > 0$ so that
\begin{equation}
\Theta_{\nu}(r) + \varepsilon < \theta.
\label{PREL31}
\end{equation}
Then, by \eqref{PREL18} and \eqref{PREL30} we get that for every $z \in D_r(0)$ there is an integer $K = K(z)$ such that
\begin{equation}
\sup_{k \geq K(z)} \left|g^{(n_k)}(z)\right|^{\frac{1}{n_k \ln n_k}} < \Theta_{\nu}(r) + \varepsilon.
\label{PREL32}
\end{equation}
It follows that if we set
\begin{equation}
G_j := \left\{z \in D_r(0) \, : \, \sup_{k \geq j} \left|g^{(n_k)}(z)\right|^{\frac{1}{n_k \ln n_k}} < \Theta_{\nu}(r) + \varepsilon\right\},
\label{PREL33}
\end{equation}
then
\begin{equation}
\bigcup_{j=1}^{\infty} G_j = D_r(0).
\label{PREL34}
\end{equation}
Thus, there is a $j_0$ for which the set $G_{j_0}$ has positive (Lebesgue) measure.

Now, as we have seen in the proof of Lemma 1, the function
\begin{equation}
\phi(z) :=  \sup_{k \geq j_0} \left|g^{(n_k)}(z)\right|^{\frac{1}{n_k \ln n_k}}
\label{PREL35}
\end{equation}
is subharmonic. Thus, as we have mentioned in Remark 3 there is an upper semicontinuous function $\hat{\phi}(z)$ such that
$\phi(z) = \hat{\phi}(z)$ for a.e. $z \in D_r(0)$. Therefore, the sets
\begin{equation*}
G_{j_0} = \left\{z \in D_r(0) \, : \, \phi(z) < \Theta_{\nu}(r) + \varepsilon\right\}
\end{equation*}
and
\begin{equation*}
\hat{G} := \left\{z \in D_r(0) \, : \, \hat{\phi}(z) < \Theta_{\nu}(r) + \varepsilon\right\}
\end{equation*}
differ by a set of measure $0$, i.e. the set $\hat{G} \vartriangle G_{j_0}$ has zero (Lebesgue) measure. Furthermore, the upper semicontinuity
of $\hat{\phi}(z)$ implies \cite{L-L} that $\hat{G}$ is open (and nonempty since $G_{j_0}$ has positive measure). Therefore, any open disk
$D_{\delta}(z_0) \subset \hat{G}$ lies almost entirely in $G_{j_0}$ in the sense that their symmetric difference has measure $0$ (in other words,
the area of $D_{\delta}(z_0) \cap G_{j_0}$ is equal to the area of $D_{\delta}(z_0)$, namely $\pi\delta^2$).

We continue by noticing that the assumption \eqref{PREL30} implies that $\theta_{\nu}(z) < \theta$ for all $z \in \mathbb{C}$ and, hence, we must have
\begin{equation}
\theta_{\mu}(z) \equiv \theta,
\label{PREL31a}
\end{equation}
where $\mu = \{m_{\ell}\}_{\ell=1}^{\infty}$ is the complementary sequence of $\nu$.

Let $D_{\delta}(z_0)$ be a disk  (with $\delta > 0$) such that $D_{\delta}(z_0) \subset \hat{G}$. If $\Gamma$ is the boundary of $D_{\delta}(z_0)$,
then due to the previous discussion we can arrange it so that the symmetric difference of $G_{j_0}$ and $\Gamma$ has one-dimensional measure $0$ (in
other words, the ``$\,$length" (i.e. the one-dimensional measure) of $\Gamma \cap G_{j_0}$ is equal to the length of $\Gamma$, namely $2\pi\delta$).

Now, by Cauchy's integral formula we have
\begin{equation}
g^{(m_{\ell})}(z_0) = \frac{1}{2\pi i} \oint_{\Gamma} \frac{g^{(n_{k[\ell]})}(z)}{(z - z_0)^{m_{\ell}-n_{k[\ell]+1}}} \, dz,
\label{PREL32a}
\end{equation}
where
\begin{equation}
k[\ell] := \max\{k \,:\, n_k < m_{\ell}\}.
\label{PREL32b}
\end{equation}
Taking absolute values in \eqref{PREL32a} yields
\begin{equation}
\left|g^{(m_{\ell})}(z_0)\right| \leq \frac{1}{2\pi \delta^{m_{\ell}-n_{k[\ell]}+1}} \oint_{\Gamma} \left|g^{(n_{k[\ell]})}(z)\right| \, ds,
\label{PREL33a}
\end{equation}
where $ds$ is the arc-length element of $\Gamma$.

By \eqref{PREL33} we have
\begin{equation}
\left|g^{(n_{k[\ell]})}(z)\right| < \left[\Theta_{\nu}(r) + \varepsilon\right]^{n_{k[\ell]} \ln n_{k[\ell]}}
\label{PREL37}
\end{equation}
for all $k \geq j_0$ and a.e. $z \in \Gamma$. Thus, by using \eqref{PREL37} in \eqref{PREL33a} we obtain
\begin{align}
\left|g^{(m_{\ell})}(z_0)\right|
&\leq \frac{1}{2\pi \delta^{m_{\ell}-n_{k[\ell]}+1}} \oint_{\Gamma} \left[\Theta_{\nu}(r) + \varepsilon\right]^{n_{k[\ell]} \ln n_{k[\ell]}} \, ds
\nonumber
\\
&= \frac{\left[\Theta_{\nu}(r) + \varepsilon\right]^{n_{k[\ell]} \ln n_{k[\ell]}}}{\delta^{m_{\ell}-n_{k[\ell]}}},
\qquad\quad
k \geq j_0,
\label{PREL38}
\end{align}
or
\begin{equation}
\left|g^{(m_{\ell})}(z_0)\right|^{\frac{1}{m_{\ell} \ln m_{\ell}}}
\leq \left(\frac{1}{\delta}\right)^{\frac{m_{\ell}-n_{k[\ell]}}{m_{\ell} \ln m_{\ell}}}
\left[\Theta_{\nu}(r) + \varepsilon\right]^{\frac{n_{k[\ell]} \ln n_{k[\ell]}}{m_{\ell} \ln m_{\ell}}},
\qquad
k \geq j_0.
\label{PREL39}
\end{equation}
Now, it is clear that
\begin{equation}
\frac{m_{\ell}-n_{k[\ell]}}{m_{\ell} \ln m_{\ell}} \to 0
\qquad \text{as }\; \ell \to \infty.
\label{B1}
\end{equation}
Also our assumption \eqref{PREL41a} for $n_k$ together with the fact that $n_{k[\ell]} < m_{\ell} < n_{k[\ell]+1}$ imply
\begin{equation}
\frac{n_{k[\ell]} \ln n_{k[\ell]}}{m_{\ell} \ln m_{\ell}} \to 1
\qquad \text{as }\; \ell \to \infty.
\label{B2}
\end{equation}
Therefore, in view of \eqref{B1} and \eqref{B2}, formula \eqref{PREL39} yields
\begin{equation}
\theta_{\mu}(z_0) = \limsup_k \left|g^{(2k+1)}(z_0)\right|^{\frac{1}{2k \ln k}} \leq \Theta_{\nu}(r) + \varepsilon < \theta,
\label{PREL40}
\end{equation}
which contradicts \eqref{PREL31a}. Hence, the assumption \eqref{PREL30} is false and we must have
$\sup_{|z| \leq r} \theta_{\nu}(z) = \theta$ for some $r > 0$, which, as we have seen earlier in the proof, implies $\theta_{\nu}(z) = \theta$
for a.e. $z \in \mathbb{C}$.
\hfill $\blacksquare$

\medskip

\textbf{Remark 4.} Suppose the complementary sequence $\mu$ of $\nu$ is also subexponential. Then, in view of \eqref{PREL18}, Theorem 1 implies
immediately that if
\begin{equation}
\mathcal{F}_{\nu} := \{z \in \mathbb{C} \,:\, \rho_{\nu}(z) = \rho\}
\qquad \text{and} \qquad
\mathcal{F}_{\mu} := \{z \in \mathbb{C} \,:\, \rho_{\mu}(z) = \rho\},
\label{PREL41}
\end{equation}
where $\rho$ is the order of $g(z)$ ($0 \leq \rho \leq \infty$) and the quantities $\rho_{\nu}(z)$ and $\rho_{\mu}(z)$ are defined by \eqref{PREL15}, then both sets $\mathcal{F}_{\nu}$ and $\mathcal{F}_{\mu}$ have full measure (and by formula \eqref{PREL17} we have
$\mathcal{F}_{\nu} \cup \mathcal{F}_{\mu} = \mathbb{C}$); in other words the sets $\mathcal{F}_{\nu}^c := \mathbb{C} \setminus \mathcal{F}_{\nu}$
and $\mathcal{F}_{\mu}^c := \mathbb{C} \setminus \mathcal{F}_{\mu}$ have Lebesque measure (i.e. area) zero.

\medskip

\textbf{Open Question.} Are the sets $\mathcal{F}_{\nu}^c$ and $\mathcal{F}_{\mu}^c$ nowhere dense in $\mathbb{C}$? Are they countable?

\section{Properties of the type}
We, now, turn our attention to the type of $g(z)$. Of course, we need to assume that $0 < \rho < \infty$.

Let $\nu = \{n_k\}_{k=1}^{\infty}$ and $\mu = \{m_{\ell}\}_{\ell=1}^{\infty}$ be two complementary sequences of indices.
Then, in view of \eqref{PREL42} we have
\begin{equation}
\tau_{\nu}(z) = \frac{1}{e \rho}\limsup_k n_k \, |a_{n_k}(z)|^{\rho/n_k}
= \frac{e^{\rho - 1}}{\rho}\limsup_k n_k^{1-\rho} \left|g^{(n_k)}(z)\right|^{\rho/n_k},
\quad
z \in \mathbb{C},
\label{PREL42t}
\end{equation}
and
\begin{equation}
\tau_{\mu}(z) := \frac{1}{e \rho}\limsup_{\ell} m_{\ell} \, |a_{m_{\ell}}(z)|^{\rho/m_{\ell}}
= \frac{e^{\rho - 1}}{\rho}\limsup_{\ell} m_{\ell}^{1-\rho} \left|g^{(m_{\ell})}(z)\right|^{\rho/m_{\ell}},
\ \
z \in \mathbb{C},
\label{PREL43}
\end{equation}
so that, as we have seen in \eqref{PREL44}, the type $\tau$ of $g(z)$ is the maximum of $\tau_{\nu}(z)$ and $\tau_{\mu}(z)$.

The following theorem gives a property of the type of $g(z)$ which is the analog of the property regarding the order of $g(z)$ established in
Theorem 1.

\medskip

\textbf{Theorem 2.} Let $\tau$ be the type of the entire function $g(z)$, while $\tau_{\nu}(z)$ is as in formula \eqref{PREL42t}, where $\nu$ is a
sequence of indices of subexponential growth.

(i) If $\tau < \infty$, then
\begin{equation}
\tau_{\nu}(z) = \tau
\qquad \text{for \;a.e. } z \in \mathbb{C}.
\label{PREL45}
\end{equation}

(ii) If $\tau = \infty$, then the set $ \{z \in \mathbb{C} \,:\, \tau_{\nu}(z) = \infty\}$ is a dense $G_{\delta}$ (therefore uncountable) subset of $\mathbb{C}$.

\smallskip

\textit{Proof}. (i) If $\tau < \infty$, we can follow the proof of Lemma 1 in order to show that $\tau_{\nu}(z)$ is subharmonic in $\mathbb{C}$. Then, by imitating the proof of Theorem 1 we can easily obtain \eqref{PREL45}.

(ii) Suppose $\tau = \infty$. Then, in view of \eqref{PREL12} we have
\begin{equation}
\sigma(z) := \sup_n n^{1-\rho} \left|g^{(n)}(z)\right|^{\rho/n} \equiv \infty.
\label{2PROOF1}
\end{equation}
Let $\mu = \{m_{\ell}\}_{\ell=1}^{\infty}$ be the complementary sequence of $\nu$. We introduce the quantities
\begin{equation}
\sigma_{\nu}(z) := \sup_k \, (n_k)^{1-\rho} \left|g^{(n_k)}(z)\right|^{\rho/n_k},
\qquad
z \in \mathbb{C}
\label{2PROOF2}
\end{equation}
and
\begin{equation}
\sigma_{\mu}(z) := \sup_{\ell} \, m_{\ell}^{1-\rho} \left|g^{(m_{\ell})}(z)\right|^{\rho/m_{\ell}},
\qquad
z \in \mathbb{C},
\label{2PROOF3}
\end{equation}
so that
\begin{equation}
\max\{\sigma_{\nu}(z), \, \sigma_{\mu}(z)\} \equiv \infty.
\label{2PROOF4}
\end{equation}
Since $\tau_{\nu}(z) = \infty$ if and only if $\sigma_{\nu}(z) = \infty$, it suffices to prove (ii) for
$\sigma_{\nu}(z)$ in place of $\tau_{\nu}(z)$.

Suppose that for some disk $D$ we had
\begin{equation}
\sup_{z \in D} \sigma_{\nu}(z) < \infty.
\label{2PROOF6}
\end{equation}
Then, \eqref{2PROOF4} would imply that $\sigma_{\mu}(z) = \infty$ for all $z \in D$, which, it can be shown to be impossible under \eqref{2PROOF6} by
following the approach used in the proof of Theorem 1, starting with formula \eqref{PREL32a}. Therefore,
\begin{equation}
\sup_{z \in D} \sigma_{\nu}(z) = \infty
\qquad \text{for any disk }\; D.
\label{2PROOF7}
\end{equation}
Now, formula \eqref{2PROOF2} implies that $\sigma_{\nu}(z)$ is lower semicontinuous on $\mathbb{C}$ (being the supremum of continuous functions).
Hence, the set
\begin{equation}
G_N := \{z \in \mathbb{C} \,:\, \sigma_{\nu}(z) > N\}.
\label{2PROOF9}
\end{equation}
is open. Furthermore, by \eqref{2PROOF7} we have that $G_N$ is dense in $\mathbb{C}$ and, therefore, the set
\begin{equation}
 \{z \in \mathbb{C} \,:\, \sigma_{\nu}(z) = \infty\} = \bigcap_{N=1}^{\infty} G_N
\label{2PROOF10}
\end{equation}
is a dense $G_{\delta}$ subset of $\mathbb{C}$.
\hfill $\blacksquare$

\medskip

\textbf{Remark 5.} Suppose that $\mu$, too, is subexponential. Then the set
\begin{equation}
\{z \in \mathbb{C} \,:\, \sigma_{\nu}(z) = \infty\} \cap \{z \in \mathbb{C} \,:\, \sigma_{\mu}(z) = \infty\},
\label{2PROOF11}
\end{equation}
being the intersection of two dense $G_{\delta}$ sets, it is again a dense $G_{\delta}$ subset of $\mathbb{C}$.

\medskip

\textbf{Remark 6.} As we have seen, the functions $\left|g^{(n)}(z)\right|^{\rho/n}$, $n = 1, 2, \ldots$, are subharmonic. It follows that $\sigma_{\nu}(z)$ satisfies \eqref{PREL21}, namely
\begin{equation}
\sigma_{\nu}(z_0) \leq \frac{1}{\pi r^2}\int_{D_r(z_0)} \sigma_{\nu}(z)\, dx dy
\label{2PROOF5}
\end{equation}
for any disk $D_r(z_0)$. Notice, however, that $\sigma_{\nu}(z)$ may not be subharmonic, since it might become infinite for some $z$ or
it might not be locally integrable.
By using \eqref{2PROOF7} in \eqref{2PROOF5}, and arguing as in the beginning of the proof of Theorem 1, we can conclude that
\begin{equation}
\int_D \sigma_{\nu}(z)\, dx dy = \infty
\qquad \text{for any disk }\; D.
\label{2PROOF8}
\end{equation}
Actually, with the help of Poisson integral formula for harmonic functions (and the fact that in any sufficiently smooth domain a subharmonic function is dominated by the harmonic function with the same boundary values) we can get a stronger version of \eqref{2PROOF8}, namely
\begin{equation}
\int_{\Gamma} \sigma_{\nu}(z)\, ds = \infty
\qquad \text{for any circle }\; \Gamma.
\label{2PROOF8a}
\end{equation}
Furthermore, since $\ln |g^{(n)}(z)|$ is subharmonic, we can work with $\ln \sigma_{\nu}(z)$ instead of $\sigma_{\nu}(z)$ and conclude that
\begin{equation}
\int_{\Gamma} \ln\sigma_{\nu}(z)\, ds = \infty
\qquad \text{for any circle }\; \Gamma.
\label{2PROOF8b}
\end{equation}
However, in spite of \eqref{2PROOF8b}, the question whether $\tau_{\nu}(z) = \infty$ for a.e. $z \in \mathbb{C}$
remains open in the case where $\tau = \infty$.

\medskip

%
%
%
%
%


\begin{thebibliography}{99}

%



%

\bibitem{H} E. Hille, \textit{Analytic Function Theory}, \textit{Volume II}, Chelsea Publishing Co., New York, N.Y., 1977.

\bibitem{L-L} E.H. Lieb and M. Loss, \textit{Analysis}, Second Edition, Graduate Studies in Mathematics, Volume 14, American Mathematical
Society, Providence, RI, 2001.


\end{thebibliography}
\end{document}